\documentclass[11pt]{amsart}

\usepackage{amssymb}
\usepackage{amsmath}
\usepackage{latexsym}
\usepackage{amsfonts}
\usepackage{graphicx}
\newtheorem{theorem}{Theorem}[section]

\newtheorem{corollary}[theorem]{Corollary}

\newtheorem{def-thm}[theorem]{Definition-Theorem}
\newtheorem{lemma}[theorem]{Lemma}

\theoremstyle{definition}
\newtheorem{definition}[theorem]{Definition}

\newtheorem{remark}[theorem]{Remark}
\newtheorem{example}[theorem]{Example}

\newtheorem{conjecture}[theorem]{Conjecture}

\newtheorem*{acknowledgement}{Acknowledgement}

\newcommand{\CC}{\mathbb{C}}

\newcommand{\rminus}{{r^-}}
\newcommand{\rmtf}{{r^{-}_{{\text{\rm tf}}}}}

\DeclareMathOperator{\kod}{kod}

\DeclareMathOperator{\Aut}{Aut}

\allowdisplaybreaks

\newcommand{\be}{\begin{equation}}
\newcommand{\ee}{\end{equation}}
\newcommand{\bea}{\begin{eqnarray}}
\newcommand{\eea}{\end{eqnarray}}

\begin{document}

\title{Reduction of manifolds with semi-negative holomorphic sectional curvature}

\author{Gordon Heier}
\address{Gordon Heier. Department of Mathematics, University of Houston, 4800 Calhoun Road, Houston, TX 77204, USA} \email{{heier@math.uh.edu}}

\author{Steven S.Y. Lu}
\address{Steven S.Y. Lu. D\'epartment de Math\'ematiques, Universit\'e du Qu\'ebec \`a Montr\'eal, C.P. 8888, Succursale Centre-Ville, Montr\'eal,  Qc, H3C 3P8, Canada} \email{{lu.steven@uqam.ca}}

\author{Bun Wong}
\address{Bun Wong. Department of Mathematics, UC Riverside, 900 University Avenue, Riverside, CA 92521, USA} \email{{wong@math.ucr.edu}}

\author{Fangyang Zheng$^1$}
\address{Fangyang Zheng. Department of Mathematics,
The Ohio State University, 231 West 18th Avenue, Columbus, OH 43210,
USA and Center for Mathematical Sciences, Zhejiang University,
Hangzhou, 310027, Zhejiang, China}
\email{{zheng.31@osu.edu}}
\thanks{$^1$Partially supported by a Simons Collaboration Grant.}

\begin{abstract}
In this note, we continue the investigation of a projective K\"ahler manifold $M$ of semi-negative holomorphic sectional curvature $H$. We introduce a new differential geometric numerical rank invariant which measures the number of linearly independent {\it truly flat} directions of $H$ in the tangent spaces. We prove that this invariant is bounded above by the nef dimension and bounded below by the numerical Kodaira dimension of $M$. We also prove a splitting theorem for $M$ in terms of the nef dimension and, under some additional hypotheses, in terms of the new rank invariant.
\end{abstract}


\maketitle

\markleft{Semi-Negative Holomorphic Sectional Curvature} \markright{Semi-Negative Holomorphic Sectional Curvature}

\section{Introduction}

In \cite{HLW_JDG}, the authors introduced a numerical invariant capturing the {\em rank} of a Hermitian manifold with semi-negative holomorphic sectional curvature. This invariant was defined as the largest codimension of maximal subspaces in the tangent spaces on which the holomorphic sectional curvature vanishes. Using this invariant, the authors established lower bounds for the nef dimension and, under certain additional assumptions, for the Kodaira dimension of projective K\"ahler manifolds. In dimension two, a precise structure theorem was obtained.\par
In this note, we continue this line of investigation, motivated by new insights about the notion of {\em rank}, which have led us to a modification of the original definition. Moreover, we will obtain new structure and splitting results as in the theorems and corollaries below. We begin by recalling the precise definition of the rank invariant from \cite{HLW_JDG}. In that paper, it was denoted by the symbol $r_M$. Here, we prefer to use the symbol $\rminus$ instead for conceptual reasons. For the sake of generality, we assume the metric to be merely Hermitian in the definition. Subsequently, we confine ourselves to the case of K\"ahler metrics.

\begin{definition}
Let $M$ be a Hermitian manifold with semi-negative holomorphic sectional curvature $H$. For $p\in M$, let $\eta(p)$ be the maximum of those integers $k\in \{0,\ldots,n:=\dim M\}$ such that there exists a $k$-dimensional subspace $L\subset T_p M$ with $H(v)=0$ for all $v\in L\backslash \{0\}$. Set $\eta_M:=\min_{p\in M} \eta(p)$ and $\rminus:=n-\eta_M$. Note that by definition $\rminus=0$ if and only if $H$ vanishes identically. Also, $\rminus=n$ if and only if there exists at least one point $p\in M$ such that $H$ is strictly negative at $p$. Moreover, $\eta(p)$ is upper-semicontinuous as a function of $p$, and consequently the set
$$\{p\in M\ |\ n-\eta(p)=\rminus\}=\{p\in M\ |\ \eta(p)=\eta_M\} $$
is an open set in $M$ (in the classical topology).
\end{definition}

The first main result in \cite{HLW_JDG} was that for a projective K\"ahler manifold of semi-negative holomorphic sectional curvature, the nef dimension satisfies $n(M) \geq \rminus$. Moreover, in \cite[Remark 1.6]{HLW_JDG} the question was raised if $n(M)$ is  necessarily equal to $\rminus$. However, the following example shows that the inequality may well be strict.

\begin{example} \label{theta_div_example}Let $A$ be a principally polarized abelian variety of dimension $n+1$, and let $M\subset A$ be a theta divisor. It was proven by Andreotti-Mayer that for generic $A$, the hypersurface $M$ is non-singular, i.e., a  submanifold. Take any flat K\"ahler metric on $A$ and restrict it to $M$. By the curvature decreasing property of subbundles, $M$ has semi-negative bisectional curvature and, in particular, $H\leq 0$. By the adjunction formula, $c_1(M)<0$, hence $n(M)=n$. On the other hand, the following explicit computation shows that $\rminus=\lfloor\frac{n+1}{2}\rfloor$, which is less than $n$ (when $n$ is at least $2$).\par
Let $p\in M$ be a fixed point. Let $U$ be a neighborhood of $p$ in $A$ with coordinates $z_0, z_1,\ldots,z_n$ such that $p$ is the origin and $T_pM=\left(\frac{\partial}{\partial z_0}\right)^\perp$. Then $z_1,\dots,z_n$ can be used as local holomorphic coordinates on $M\cap U$ and $M$ defined as a graph by $z_0=f(z_1,\ldots,z_n)$ for some holomorphic function $f$ with $f(0,\ldots,0)=0$ and $df(0,\ldots,0)=0$. The induced metric on $M$ is the graph metric which is given in components as
$$g_{i\overline{j}} = \delta_{i j} +f_i\overline{f_j},$$
where $f_i = \frac{\partial f}{\partial z_i}$. From this, we obtain that the components of the curvature tensor at the point $p$ are
\begin{equation}\label{theta_div_curv_tensor} R_{i\overline{j} k \overline{l}}=-f_{ik}\overline{f_{jl}},\end{equation}
where $f_{ij} = \frac{\partial^2 f}{\partial z_i\partial z_j}$. Therefore, for a unit vector $v=
\sum_{i=1}^n v_i \frac{\partial}{\partial z_i} \in T_pM$, the holomorphic sectional curvature of $v$ is
$$R_{v\overline{v}v\overline{v}}=-|f_{vv}|^2,$$
where $f_{vv} = \sum_{i,j=1}^n v_iv_j f_{ij}$. Thus, for any non-zero vector $v\in T_pM$, $H(v)=0$ if and only if $f_{vv}=0$. Note that all such $v$ form a quadratic cone $Q$ in $T_p M=\CC^n$. According to \cite[page 735, Proposition]{GH}, every quadratic cone of dimension $m$ contains a linear space of dimension $\lfloor\frac{m-1}{2}\rfloor+1$, which in our situation yields $\rminus\leq \lfloor\frac{n+1}{2}\rfloor$ (note that $m=n-1$). Moreover, if the quadratic cone is non-degenerate, there exist no linear spaces of larger dimension on it. We will conclude the discussion of this example by showing that there exists a point $p$ of $M$ where $Q$ is indeed non-degenerate and can therefore infer that $\rminus= \lfloor\frac{n+1}{2}\rfloor$.\par
Let us assume that the matrix $(f_{ij})$ is degenerate at $p$. Then there exists a non-zero tangent vector $v\in T_pM$ such that $f_{iv}:=\sum_{j=1}^n v_j f_{ij}=0$ for all $i=1,\ldots,n$. Consequently, by \eqref{theta_div_curv_tensor}, $R_{i\overline{j}v\overline{v}}=0$ for all $i,j=1,\ldots,n$, which implies that the Ricci curvature form evaluated on $v$ is zero. In particular, the $n$-fold wedge of the Ricci curvature form with itself vanishes at $p$. Now, if this holds for all $p\in M$, then $c_1(M)^n =0$, contradicting the initial observation that $c_1(M)<0$. \qed

\end{example}

The problem here is that, for a compact K\"ahler manifold $M$ with $H\leq 0$, there could be many zeros of $H$, even when $M$ has ample canonical line bundle (and even when the bisectional curvature of $M$ is semi-negative). These zeros, as in the above example of the theta divisor, may actually not prevent the Ricci curvature from being negative and thus do not contribute towards the defect in nef dimension (or Kodaira dimension). In order to capture the differential geometric aspect of this phenomenon, let us introduce the following definition.

\begin{definition} Let $M$ be a K\"ahler manifold with semi-negative holomorphic sectional curvature, and $p\in M$. A (type $(1,0)$ complex tangent) vector $v\in T_pM$ is called {\em truly flat} if  for any $w\in T_pM$ that is perpendicular to $v$, it holds $R_{u\overline{u}u\overline{u}}=R_{w\overline{w}w\overline{w}}$, where $u=v+w$.
\end{definition}

Under the condition $H\leq 0$, a polarization argument similar to the one in \cite[Section 2]{WZ_Contemp} proves that $v\in T_pM$ is truly flat if and only if $R_{v\overline{x}y\overline{z}}=0$ for any $x$, $y$, and $z$ in $T_pM$, in other words, if and only if $v$ belongs to the kernel of $R$ (aka nullity of the curvature).\par

\begin{definition}
Let $M$ be a K\"ahler manifold with semi-negative holomorphic sectional curvature.  Let $V_p\subset T_pM$ be the set of all truly flat vectors at $p$. It is a complex linear subspace. Define $\rmtf$ (also written as $\rmtf(M)$ in case of ambiguity) to be the maximal codimension of $V_p$ over all $p\in M$. 
\end{definition}

Clearly, our two rank invariants satisfy $\rminus\leq \rmtf$. In Example \ref{theta_div_example}, the ampleness of the canonical line bundle of $M$ implies $n=\rmtf$ and thus $$\rminus=\lfloor\frac{n+1}{2}\rfloor < n = \rmtf$$ (for $n\geq 2$). Our first main result in this note is that Theorem 1.4 of \cite{HLW_JDG} actually holds for the (possibly larger) invariant $\rmtf$.

\begin{theorem} \label{upper_bound_thm} Let $M$ be a projective manifold with a K\"ahler metric whose holomorphic sectional curvature is semi-negative. Then its nef dimension $n(M)$ satisfies
$$ n(M) \geq \rmtf.$$ \end{theorem}

Our next result is the following splitting theorem, which is similar in spirit to the results in \cite{WZ_JDG} and \cite{Liu}. Its proof first establishes the fact that under the assumption of semi-negative holomorphic sectional curvature, the nef fibration is actually a metric holomorphic fiber bundle up to an unramified covering. The actual product decomposition is then obtained by an application of Serre's discussion \cite{Serre} of {\it espaces fibr\'es principaux} in the projective setting. It should be noted that, in contrast to the results in \cite{WZ_JDG} and \cite{Liu}, our theorem does not require an assumption of semi-negativity for the {\it bisectional} curvature.

\begin{theorem} \label{split_thm} Let $M$ be a projective manifold with a K\"ahler metric whose holomorphic sectional curvature is semi-negative. Then there exist a finite unramified cover $\widehat M$ of $M$, a flat abelian variety $A$ and  a projective K\"ahler manifold $N$ of semi-negative holomorphic sectional curvature with $\dim N = n(N)=n(M)$, such that $\widehat M$ is biholomorphic and isometric to $A\times N$ and locally isometric to $M$.
\end{theorem} 

As a counterpoint to Theorem \ref{upper_bound_thm}, the following theorem establishes that $\rmtf$ is bounded below by the numerical Kodaira dimension $\nu (M)$. 

\begin{theorem} \label{lower_bound_thm} Let $M$ be a projective manifold with a K\"ahler metric whose holomorphic sectional curvature is semi-negative. Then $$\rmtf \geq \nu (M).$$ \end{theorem}

Theorems \ref{upper_bound_thm} and \ref{lower_bound_thm} allow us to predict that $\rmtf$ is the correct invariant to consider, since $n(M)$ and $\nu(M)$ are commonly believed to always agree. In fact, this agreement is a special case of the (much stronger) Abundance Conjecture which states that the known chain of inequalities $$n(M) \geq \nu(M)\geq \kod(M)$$ is actually a chain of equalities. The Abundance Conjecture is known up to dimension $3$ by the works of Miyaoka and Kawamata. We refer the reader to \cite[Section 4]{HLW_JDG} and \cite[Part I, Lecture IV]{Miyaoka_Peternell_book} for more information in this direction. We may sum up Theorems \ref{upper_bound_thm} and \ref{lower_bound_thm} in the form of the following succinct immediate corollary.
\begin{corollary}\label{chain_inequ}
Let $M$ be a projective manifold with a K\"ahler metric whose holomorphic sectional curvature is semi-negative. Then
$$n(M) \geq \rmtf \geq \nu(M)\geq \kod(M).$$
If the Abundance Conjecture holds true for $M$, then the line above is an all-around equality.
\end{corollary}

The effect of the Abundance Conjecture on Theorem \ref{split_thm} is expressed in the following corollary.

\begin{corollary} \label{split_thm_ac} Let $M$ be a projective manifold with a K\"ahler metric whose holomorphic sectional curvature is semi-negative. Assume that the Abundance Conjecture holds in dimension $n(M)$. Then there exist a finite unramified cover $\widehat M$ of $M$, a flat abelian variety $A$ and  a projective K\"ahler manifold $N$ of semi-negative holomorphic sectional curvature with $c_1(N) < 0$ and $\dim N =\rmtf(M)$, such that $\widehat M$ is biholomorphic and isometric to $A\times N$ and locally isometric to $M$.
\end{corollary} 

The proof of this corollary merely consists of a few standard steps in addition to the proof of Theorem \ref{split_thm}. Namely, Abundance yields $\dim N=n(N)=\kod(N)$, i.e., $N$ is of general type. Since $N$ retains the semi-negative holomorphic sectional curvature, it does not contain rational curves (see \cite[Lemma 3.1]{HLW_JDG} and the comments preceding it). It is then well-known (see \cite[p.\ 219, Exercise 8]{Debarre}) that, due to work of Kawamata, the canonical line bundle of $N$ is ample, i.e., $c_1(N) < 0$. Finally, $c_1(N) < 0$ implies that $\rmtf(M)= \dim N$. \par
We would like to emphasize that Corollary \ref{split_thm_ac} generalizes \cite[Theorem 1.10]{HLW_JDG} from the case of  surfaces to threefolds.

Under the additional assumption of the bisectional curvature being semi-negative, we note that the conclusion of Corollary \ref{split_thm_ac} holds without any assumption of Abundance. In fact, analogous statements were proven in \cite[Theorem E]{WZ_JDG} (under the additional assumption that the metric is real-analytic) and \cite{Liu} in terms of the Ricci rank, and their proofs go through also in our case, yielding the following theorem. Recall that the holomorphic sectional curvature is simply the restriction of the bisectional curvature to the diagonal, so compared to Theorem \ref{split_thm}, we have to work with a strictly stronger curvature assumption in this theorem to be able to obtain the result. 

\begin{theorem} \label{split_thm_ii} Let $M$ be a projective manifold with a K\"ahler metric whose bisectional curvature is semi-negative. Then there exist a finite unramified cover $\widehat M$ of $M$, a flat abelian variety $A$ and  a projective K\"ahler manifold $N$ of semi-negative holomorphic sectional curvature with $c_1(N) < 0$ and $\dim N =\rmtf(M)$, such that $\widehat M$ is biholomorphic and isometric to $A\times N$ and locally isometric to $M$.

\end{theorem} 

Finally, we would like to remark that there has recently been strong renewed interest and very significant progress in the study of holomorphic sectional curvature, both of (semi-)negative and (semi-)positive sign. Most relevantly to the present article, there has appeared the sequence of works \cite{HLW_MRL}, \cite{Wu_Yau_Invent}, \cite{Tosatti_Yang}, \cite{Wu_Yau_CAG}, and \cite{Diverio_Trapani}. The first two papers focussed on Yau's original conjecture \cite[Conjecture 1]{Wu_Yau_Invent} stating that a projective K\"ahler manifold with negative holomorphic sectional curvature has ample canonical line bundle, with \cite{HLW_MRL} proving the case of the dimension being no greater than 3 and \cite{Wu_Yau_Invent} proving the general case. In the paper \cite{Tosatti_Yang}, the authors were able to establish the conjecture with the projectivity assumption removed. The subsequent papers \cite{Wu_Yau_CAG} and \cite{Diverio_Trapani} then established the conjecture under the further weakened assumption of quasi-negativity of the holomorphic sectional curvature. In our language, this is equivalent to the assumption of semi-negative holomorphic sectional curvature and $r^-=\dim M$.\par
In light of Example \ref{theta_div_example} and the introduction and discussion of $\rmtf$ in this paper, the authors believe that the following conjecture represents an optimal extension of Yau's original conjecture. Note that the conjecture is essentially dictated by Corollary \ref{chain_inequ} and is thus a theorem in dimension no greater than 3 in the projective case. Due to Theorem \ref{split_thm_ii}, we know the conjecture to be true in the projective case under the assumption of semi-negative bisectional curvature. Note also that an earlier version of this conjecture appeared as \cite[Conjecture 4.3]{HLW_JDG}.
\begin{conjecture}
Let $M$ be a compact K\"ahler manifold with semi-negative holomorphic sectional curvature. Then the equality $\kod(M) = \rmtf$ holds. In particular, if $\rmtf = \dim M$, then the canonical line bundle of $M$ is ample. If $M$ is additionally assumed to be projective, then the equality $\kod(M) = \rmtf = n(M)$ holds, i.e., Abundance holds on $M$. \end{conjecture}

\section{The proofs of the Theorems}

In this section, we will provide the proofs of the theorems stated in the introduction. We will use the setup and notations as in \cite{HLW_JDG}.\par

\noindent {\bf Proof of Theorem \ref{upper_bound_thm}.}  Let $M$ be an $n$-dimensional projective manifold with a K\"ahler metric whose holomorphic sectional curvature $H$ is semi-negative, i.e., $H\leq 0$ everywhere. Assume that $H$ is not identically zero (otherwise $M$ is a flat complex space form). By Theorem 1.4 in \cite{HLW_JDG}, the canonical line bundle $K_M$ is nef, with positive numerical dimension, and the nef reduction map $f: M \dashrightarrow Y$ will have the generic compact fibers being flat abelian varieties up to an unramified covering, of dimension $k=n-n(M)$.\par

Let $A$ be a generic fiber, which up to an unramified covering is a flat abelian variety of dimension $k$. Fix $p\in A$. For any $v\in T_pA$, we have $R_{v\overline{v}v\overline{v}}=0$. Here and below, $R$ stands for the curvature tensor of $M$. By Part (1) of Lemma \ref{key_lemma} below, we know that $R_{v\overline{v} w\overline{w}} \leq 0$ for any $w\in T_pM$. So the Ricci tensor $R_{v\overline{v}}\leq 0$ along $A$. Since the Ricci $(1,1)$-form represents the anti-canonical line bundle, while $K_M|_A = K_A$ is trivial, we conclude that $R_{v\overline{v}w\overline{w}} = 0$ for any $w\in T_pM$. By Part (2) of Lemma \ref{key_lemma} below, this implies that $R_{v\overline{x}y\overline{z}}=0$ for any $x$, $y$, and $z$ in $T_pM$. Therefore, $v$ represents a truly flat direction of $H$. So $T_pA\subset V_p$, thus $k \leq n-\rmtf$ and $n(M) \geq \rmtf$. This completes the proof of Theorem \ref{upper_bound_thm}. \qed

\begin{lemma} \label{key_lemma} Let $M$ be a K\"ahler manifold with semi-negative holomorphic sectional curvature $H$, and let $p\in M$, $0\neq v\in T_pM$.

(1). If $H(v)=0$, then $R_{v\overline{v} w\overline{w}} \leq 0$ for any $w\in T_pM$.

(2). If $R_{v\overline{v} w\overline{w}} = 0$ for any $w\in T_pM$, then $R_{v\overline{x}y\overline{z}}=0$ for any $x$, $y$, and $z$ in $T_pM$.
\end{lemma}
\proof Consider the quantity $F(\lambda ):= R_{u\overline{u} u\overline{u}}\leq 0$ with $u=\lambda v +w$, where $\lambda $ is a complex number. Since $R_{v\overline{v} v\overline{v}} = 0$,  $F(\lambda )$ is equal to
$$ 2Re(|\lambda |^2\lambda R_{v\overline{v} v\overline{w}}) + 4|\lambda |^2 R_{v\overline{v} w\overline{w}} + 2Re(\lambda^2 R_{v\overline{w} v\overline{w}}) + 2Re(\lambda R_{v\overline{w} w\overline{w}}) + R_{w\overline{w} w\overline{w}}.$$
So the non-positivity of $F$ implies that
$R_{v\overline{v} v\overline{w}}=0$ and
$$2R_{v\overline{v} w\overline{w}} + |R_{v\overline{w} v\overline{w}}| \leq 0, \ \ \ R_{w\overline{w} w\overline{w}}\cdot (2R_{v\overline{v} w\overline{w}} + |R_{v\overline{w} v\overline{w}}|) \geq |R_{v\overline{w} w\overline{w}}|^2.$$
In particular, $R_{v\overline{v} w\overline{w}} \leq 0$ so part (1) holds. Also, when $R_{v\overline{v} w\overline{w}} =0$, the above inequalities lead to $R_{v\overline{w} w\overline{w}} =0$. By polarization, one gets $R_{v\overline{x}y\overline{z}}=0$ for any $x$, $y$, $z$ in $T_pM$, so part (2) holds as well. \qed

Next, let us prove Theorem \ref{split_thm}, following the philosophy of \cite{WZ_JDG}. Note again that in that paper, the authors worked under the additional assumption that the bisectional curvature was semi-negative. Here, we are in the better situation that the generic fibers of the nef reduction map are compact and can thus proceed without this assumption.

\noindent {\bf Proof of Theorem \ref{split_thm}.} Let $M$ be an $n$-dimensional projective manifold equipped with a K\"ahler metric with semi-negative holomorphic sectional curvature $H$. Let $n(M)$ be the nef dimension of $M$, and $f: M \dashrightarrow Y$ be the nef reduction map. Let $F=F_p$ be the fiber of $f$ at a generic point $p\in M$. By Theorem 1.4 of \cite{HLW_JDG}, we know that $F$ is the finite unramified undercover of a $(n-k)$-dimensional abelian variety, where $k=n(M)$, and the restriction metric on $F$ is flat. In fact, for the nef reduction map $f$, there are Zariski open subsets $M^{\ast} \subset M$ and $Y^{\ast }\subset Y$ such that $h=f|_{M^{\ast}}: M^{\ast} \rightarrow Y^{\ast}$ is a holomorphic submersion with compact fibers. \par

First let us show that $M^{\ast }$ is locally a product space, namely, given any $p\in M^{\ast }$, there exists an open neighborhood $V$ of $h(p)\in Y^{\ast }$ such that $h^{-1}(V)$ is holomorphically isometric to the product $F\times V$. In particular, all fibers of $h$ are isomorphic to each other. \par

To see this, let us choose a small coordinate neighborhood $V$ of $h(p)$ in $Y^{\ast }$ so that $V$ is biholomorphic to a polydisc and $U=h^{-1}(V)$ is diffeomorphic to the product space $F\times V$ by Ehresmann's theorem. The universal cover $\widetilde{U}$ of $U$ is a holomorphic fibration over $V$ with each fiber being ${\mathbb C}^{n-k}$. So $\widetilde{U}$ is a holomorphic fiber bundle over $V$, and replacing $V$ by a smaller neighborhood if necessary, we may assume that $\widetilde{U}$ is biholomorphic to $ {\mathbb C}^{n-k}\times V$. Let $z=(z_1, \ldots z_k)$ be a local holomorphic coordinate in $V$ centered at $h(p)$, and let $t=(t_1, \ldots , t_{n-k})$ be the standard Euclidean coordinate in ${\mathbb C}^{n-k}$. Then $(t,z)$ gives a holomorphic coordinate on $\widetilde{U}= {\mathbb C}^{n-k}\times V$. \par

Denote by $\Gamma$ the deck transformation group with $U=\widetilde{U}/\Gamma$. There is a normal subgroup of $\Gamma '\subseteq \Gamma $ of finite index such that $\Gamma ' \cong {\mathbb Z}^{2n-2k}$, and each element in it is in the form $\gamma (t,z) = (t+v_{\gamma }(z),z)$. For any fixed $z$, the vectors $v_{\gamma }(z)$ for all $\gamma \in \Gamma '$ form a lattice in ${\mathbb C}^{n-k}$, and each $v_{\gamma }(z)$ depends holomorphically on $z$ as the fibers of the nef reduction map $h$ vary holomorphically. Let us fix a set of generators $\{ \gamma_1, \ldots , \gamma_{2n-2k} \}$ for $\Gamma'$, and assume that $\{ v^{1}(z), \ldots , v^{{n-k}}(z) \} $ is linearly independent over ${\mathbb C}$ for each $z\in V$. Here, we wrote $v^{\alpha }$ for $v_{\gamma_{\alpha }}$. For $1\leq \alpha \leq n-k$, write $e_{\alpha} = (0, \ldots , 1, \ldots , 0) \in {\mathbb C}^{n-k}$, where the only non-zero entry is $1$ and is at the $\alpha$-th position. Denote by $P$ the $(n-k)\times (n-k)$ matrix of holomorphic functions whose  $\alpha$-th row is given by $v_{{\alpha }}(z)$ for each $1\leq \alpha \leq n-k$. Consequently, we have $e_{\alpha }P = v^{{\alpha}}(z)$ for each $\alpha$. \par

So replacing the coordinates $(t,z)$ by $(\tilde{t}(t,z),z)$ with $\tilde{t}=tP^{-1}$ if necessary, we may assume that $v^{\alpha}(z)=e_{\alpha}$ for the first $n-k$ elements of the generator set $\{ \gamma_1, \ldots , \gamma_{2n-2k} \}$ of $\Gamma'$. The $v$ for the next $n-k$ elements will be denoted by $w_1(z), \ldots , w_{n-k}(z)$. For any fixed $z$ in $V$, the set $\{ e_1, \ldots , e_{n-k}, w_1, \ldots , w_{n-k}\}$ of vectors in ${\mathbb C}^{n-k}$ is linearly independent over ${\mathbb R}$. \par

Let us express the pull back of the K\"ahler metric $g$ on $\widetilde{U}$ in terms of the coordinates $(t,z)$:

\begin{align*} \omega_g\ =\ & \sqrt{-1}\Big( \sum_{\alpha , \beta } g_{\alpha \overline{\beta }} dt_{\alpha }\wedge d\overline{t}_{\beta } +\sum_{i,\alpha }g_{n-k+i\, \overline{\alpha }} dz_i\wedge d\overline{t}_{\alpha } \\ & + g_{\alpha\, \overline{n-k+i}} dt_{\alpha }\wedge d\overline{z}_i  +\sum_{i,j} g_{n-k+i\, \overline{n-k+j}} dz_i\wedge d\overline{z}_j \Big).
\end{align*}

Here, we adopted the index range convention $1\leq i,j \leq k$ and  $1\leq \alpha , \beta \leq n-k$. Since $\nabla _{\partial_{\alpha }}\partial_{\beta} =0$, and the connection matrix under the natural frame is given by $\partial g g^{-1}$, we know that $\partial_{\alpha }g_{\beta \overline{\ast }}=0$. So each $g_{\alpha \overline{\beta }}=g_{\alpha \overline{\beta }} (z,\overline{z})$ is independent of $t$, and each $g_{\alpha\, \overline{n-k+i}}$ is anti-holomorphic in $t$. By the K\"ahlerness of the metric, we know that the partial derivatives satisfy
$  g_{\alpha\, \overline{n-k+i}, \overline{\beta} } =  g_{\alpha \overline{\beta }, \overline{n-k+i}} $, 
which is independent of $t$, so we know that
$$ g_{\alpha\, \overline{n-k+i}} = \sum_{\beta}  g_{\alpha \overline{\beta}, \overline{n-k+i}} \ \overline{t}_{\beta} + f_{\alpha\,\overline{n-k+i}}(z,\overline{z}). $$
Similarly, we get from the K\"ahlerness of $g$ the following:

\begin{align*}
g_{n-k+i\, \overline{n-k+j}}\ =\ & \sum_{\alpha , \beta } g_{\alpha \overline{\beta } ,n-k+i\, \overline{n-k+j}} \ t_{\alpha } \overline{t}_{\beta } + \sum_{\alpha } f_{\alpha\, \overline{n-k+j},n-k+i} \ t_{\alpha }\\& + \sum_{\alpha} f_{n-k+i\, \overline{\alpha},\overline{n-k+j}} \ \overline{t}_{\alpha} +  \xi_{i\overline{j}} (z,\overline{z}).
\end{align*}

For any deck transformation $\gamma (t,z)=(t+v(z),z)$, we have  $\gamma^{\ast }\omega_g = \omega_g$. By comparing the coefficients in front of $dt_{\alpha}\wedge d\overline{z}_i$, we get $ g_{\alpha\, \overline{n-k+i}} = \gamma^{\ast } g_{\alpha\, \overline{n-k+i}} + \gamma^{\ast } g_{\alpha \overline{\beta}} \  \overline{v_{\beta,n-k+i}} $, or equivalently, $ (g_{\alpha \overline{\beta}} \ \overline{v_{\beta}})_{\overline{n-k+i}} =0$. Now if we let $\gamma$ be $\gamma_{\beta}$, our $v$ is $e_{\beta}$, so $v_{\beta_1}=\delta_{\beta_1 \beta}$, thus we get $g_{\alpha \overline{\beta }, \overline{n-k+i}}=0$ for any $\alpha$, $\beta$ and any $i$. So each $g_{\alpha \overline{\beta}}$ is a constant. \par

So now we know that each $g_{\alpha \overline{\beta}}$ is a constant, each $g_{\alpha\, \overline{n-k+i}}$ is a function of $z$  and $\overline{z}$, and
$$ g_{n-k+i\, \overline{n-k+j}} = \sum _{\alpha } g_{\alpha\, \overline{n-k+j}, n-k+i} \ t_{\alpha } +  \sum_{\beta } g_{n-k+i\, \overline{\beta}, \overline{n-k+j}} \ \overline{t}_{\beta} + \ \xi_{i\overline{j}}(z,\overline{z}).$$
Similarly, by comparing the term $dz_i\wedge d\overline{z}_j$ in $\gamma^{\ast }\omega_g = \omega_g$, we get
$$(g_{\alpha\, \overline{n-k+j}}\ v_{\alpha})_{,n-k+i} + (g_{n-k+i\, \overline{\alpha}}\  \overline{v}_{\alpha})_{,\overline{n-k+j}} =0.$$ That is, we have 
$$g_{v\, \overline{n-k+j},n-k+i}+g_{n-k+i\, \overline{v},\overline{n-k+j}}=0,$$ 
or equivalently, $(g_{n-k+i\, \overline{n-k+j}})_{,v+\overline{v}}=0 $
for the vector $v$ associated with any deck transformation $\gamma$ in $\Gamma'$. By taking $\gamma$ from the aforementioned generator set $\{\gamma_1, \ldots , \gamma_{2n-2k}\}$, as $\{e_1, \ldots , e_{n-k},\allowbreak w_1, \ldots , w_{n-k}\}$ is linearly independent over ${\mathbb R}$, we know that the square matrix formed by the imaginary parts of $w_i$ must be non-singular, thus the above condition yields 
$$g_{n-k+i\, \overline{n-k+j}, \alpha}=0.$$ 
So the entire matrix of the metric $g$ is independent of $t$. \par

Let us denote by  ${\mathcal F}$ the totally geodesic foliation in $M^{\ast }$ given by the fibers of the nef reduction map $h$, and denote by ${\mathcal F}^{\perp}$ its orthogonal complement distribution. Then by the fact that each $g_{\alpha \overline{\beta}}$ is a constant, we know that the twisting tensor of ${\mathcal F}$ vanishes, so ${\mathcal F}^{\perp}$ is also a totally geodesic foliation. So the metric $g$ splits as a product metric locally in $M^{\ast }$. In particular, the fibers of $h$ vary in a parallel way, thus they are all isomorphic to each other since  $M^{\ast }$ is path connected. \par

Now since $M$ is the closure of $M^{\ast }$, for any $p\in M\setminus M^{\ast }$, there is a sequence of points $\{ p_i\}$ in $M^{\ast }$ approaching $p$. The limiting positions of $F_{p_i}$ must be unique since $F_q$ are parallel for all $q \in M^{\ast }$. This means that $M$ has a flat de Rham factor: the universal cover $\widetilde{M}$ is holomorphically isometric to ${\mathbb C}^{n-k}\times X^k$, where $k=n(M)$ and $X$ is a complete K\"ahler manifold of dimension $k$ with semi-negative holomorphic sectional curvature. \par

As in the proof of Theorem E in \cite{WZ_JDG}, let us denote by $\Gamma $ the deck transformation group of $\widetilde{M}$, and by $I_1$, $I_2$ the groups of holomorphic isometries of ${\mathbb C}^{n-k}$ and $X$, respectively. Since $\Gamma$ needs to preserve the vertical foliation of $f$ and its perpendicular complement,  each element $\gamma$ of $\Gamma$ is in the form $(\gamma_1, \gamma_2)$, where $\gamma_i\in I_i$. Denote by $\Gamma_i$ the image of $\Gamma$ under the projection of $I_1\times I_2$ onto its factors. The leaves ${\mathbb C}^{n-k}\times \{ x\}$ close up in $M$, so the corresponding group $\Gamma_2$ must be discrete. Passing to a finite index subgroup of $\Gamma$, we may assume that $\Gamma_2$ acts freely on $X$, thus giving a projective K\"ahler manifold $B=X/\Gamma_2$, and we get  the holomorphic fiber bundle $f_1: M' \rightarrow B$ with the fiber $A$ being an abelian variety and $M'$ being a finite unramified cover of $M$. Moreover, $f_1$ is a metric bundle, i.e., any point in $B$ is contained in a neighborhood $U$ such that $f_1^{-1}(U)$ is isometric to $A\times U$.\par

In fact, a stronger statement is true in our projective setting. Namely, for any given point $b\in B$, we can take hyperplane sections of $M'$ to find a Zariski neighborhood $U(b)$ in the base and a multisection of $f_1$ over $U(b)$ such that the restriction of the bundle projection over $U(b)$ to the multisection is finite \'etale. It follows that the fibration given by $f_1$ over $U(b)$ becomes trivial after the base change by this finite \'etale map.  Therefore, the map $f_1$ is {\it locally isotrivial} and thus an {\it espace fibr\'e principal} in the sense of \cite[2.2 D\'efinitions]{Serre}, which is also known as a {\it torsor} (\cite[p.\ 120]{Milne}). Now, the isomorphism classes of torsors over $B$ with fiber $A$ are parameterized by a certain first cohomology group of $B$ with coefficients in $A$ (\cite[Proposition 1]{Serre}), which is a torsion group in our smooth projective setting (\cite[Lemme 7]{Serre}). A standard argument (\cite[Proposition 17]{Serre}) then shows, by a base change $N\to B$ via the ``$m$-torsion points of the fibers" where $m$ is the order of our torsor, that there exists an unramified cover $\widehat M= M'\times_B N$ of $M'$ such that $\widehat M$ is biholomorphic to the product $A\times N$ and necessarily isometric, with $N$ being a projective K\"ahler manifold of semi-negative holomorphic sectional curvature.  \par
Since nef dimension is clearly invariant under taking unramified covers and additive in products, we have 
\begin{equation*}
\pushQED{\qed} 
n(N)=n(B)=n(\widehat M)=n(M)=\dim N. \qedhere
\popQED
\end{equation*}
~\\[-.8cm]

\begin{remark}
The proof of Theorem \ref{split_thm} yields a non-existence result for K\"ahler metrics of semi-negative holomorphic sectional curvature on the total space of families of abelian varieties over a projective base of general type whose fibers are not all isomorphic to each other. To be more precise, let $M$ be a projective manifold and let $f:M\to B$ be a holomorphic submersion over a projective base manifold $B$ of general type such that the fibers are abelian varieties which are not all isomorphic to each other. According to the definition, $f$ is a nef reduction map of $M$. If $M$ were to admit a K\"ahler metric $H$ of semi-negative holomorphic sectional curvature, then the restriction of $H$ to any fiber of $f$ is necessarily flat and totally geodesic in $M$. Moreover, the proof of Theorem \ref{split_thm} then shows that all fibers are holomorphically isometric to each other, and $f$ is necessarily a holomorphic fiber bundle. This is a contradiction.
\end{remark}

\noindent {\bf Proof of Theorem \ref{lower_bound_thm}.} Let us write $k$ for $\rmtf$ and set $n=\dim M$. Then at any point $p\in M$, there are at least $n-k$ linearly independent truly flat tangent vectors with respect to $H$. By Lemma \ref{key_lemma}, we know that there is a unitary local tangent frame $\{ e_1, \ldots  , e_n\}$ near $p$ such that $R_{i\overline{\alpha}\beta\overline{\gamma}}=0$ for any $k+1\leq i\leq n$ and any $1\leq \alpha , \beta , \gamma \leq n$. So if we denote by $\{ \varphi_1, \ldots , \varphi_n\}$ the coframe dual to $e$, then the Ricci $(1,1)$-form of $M$ will be
$$ \xi = \sum_{i,j=1}^k  R_{i\overline{j}} \varphi_i \wedge \overline{\varphi}_j ,$$
where $R_{i\overline{j}} = \sum_{\alpha =1}^n R_{i\overline{j}\alpha\overline{\alpha}}$ are the components of the Ricci curvature tensor.

From this, we see that the wedge product of $(k+1)$ copies of $\xi$ is identically zero on $M$, i.e., $\xi^{k+1}=0$.  Since $\xi$ represents $c_1(M)$, this implies that $c_1(M)^{k+1}=0$.  Therefore, the numerical Kodaira dimension $\nu(M)$ is no greater than $k$. This completes the proof of Theorem \ref{lower_bound_thm}.
 \qed

\noindent {\bf Proof of Theorem \ref{split_thm_ii}.} In the present situation, the polarization argument from \cite[Section 2]{WZ_Contemp} proves that $v\in T_pM$ is truly flat if and only if $v$ is in the kernel of the Ricci form if and only if $R_{v\overline{x}y\overline{z}}=0$ for any $x$, $y$, and $z$ in $T_pM$. Therefore, the proof of Theorem 2 in \cite{Liu} applies to our case as well. \par

We take this opportunity to make a very minor correction to the penultimate paragraph on \cite[p.~280]{WZ_JDG}. Let $U$ be the real Lie algebra of $A^*\subset \Aut(\widetilde M)$. Then $U$ is a real vector space with a basis $u_1,\ldots,u_k$. Each vector $u_i$ is 
the real part of a holomorphic vector field $v_i:= u_i-\sqrt{-1}Ju_i$, $i=1,\ldots,k$. Then the vector space $V$ spanned by these $v_1\ldots,v_k$ is the $V$ used in \cite{WZ_JDG}.\qed

We conclude by returning to the topic of Abundance and Kodaira dimension with the following remark.
\begin{remark}
Let us denote by ${\mathcal H} = {\mathcal H}^{0-}_n$ the set of all $n$-dimensional projective manifolds $M$ that are equipped with K\"ahler metrics of semi-negative holomorphic sectional curvature $H$. For such a manifold, its canonical line bundle is nef, and Corollary \ref{chain_inequ} established the chain of inequalities $\kod(M)\leq \nu(M)\leq \rmtf \leq n(M)$. This suggests that, perhaps ${\mathcal H}$ would be a good test case for the  Abundance Conjecture, since if it fails on such a manifold, the curvature rank $\rmtf$ would either be larger than $\kod(M)$ or smaller than $n(M)$, and either way, one has some additional differential geometric tools to work with. Moreover, since the condition $H\leq 0$ naturally rules out the presence of rational curves, which is at the center of all birational tanglement, we feel that this class should be the first trial case for pushing the deformation theorem of Huai-Dong Cao on K\"ahler-Ricci flow from the $c_1<0$ case to a more general case, in which $\kod(M)$ may be less than $n$.
\end{remark}

\begin{acknowledgement}
We would like to thank Kefeng Liu, Hongwei Xu, and Shing-Tung Yau for their interest.
\end{acknowledgement}

\end{document}